# Product and anti-Hermitian structures on the tangent space

E. Peyghan, A. Razavi and A. Heydari
August 22, 2007

## Abstract

Noting that the complete lift of a Rimannian metric $g$ defined on a differentiable manifold $M$ is not 0-homogeneous on the fibers of the tangent bundle $TM$. In this paper we introduce a new lift $\tilde{g}_2$ which is 0-homogeneous. It determines on slit tangent bundle a pseudo-Riemannian metric, which depends only on the metric $g$. We study some of the geometrical properties of this pseudo-Riemannian space and define the natural almost complex structure and natural almost product structure which preserve the property of homogeneity and find some new results.

## 1. Introduction.

The importance of the complete lift $g_2$ of a Riemannian metric $g$ is well known in Riemannian geometry, Finsler geometry and Physics, and has many applications in Biology too (see [1]). The tensor field $g_2$ determines a pseudo-Riemannian structure on slit tangent bundle $\widetilde{TM} = TM \setminus \{0\}$, but $g_2$ is not 0-homogeneous on the fibers of the tangent bundle $TM$. Therefore, we cannot study some global properties of the pseudo-Riemannian space $(\widetilde{TM}, g_2)$. For instance we can not prove a theorem of Gauss-Bonnet type for this space (see [4]).

In this paper, we define a new kind of lift $\tilde{g}_2$ to $TM$ of the Riemannian metric $g$. Thus $\tilde{g}_2$ determines on $\widetilde{TM}$ a pseudo-Riemannian structure, which is 0-homogeneous on the fibers of $TM$ and depends only on $g$. Some geometrical properties of $\tilde{g}_2$ such as the Levi-Civita connection are studied.

Almost complex and almost product structures are among the most important geometrical structures which can be considered on a manifold. Geometric properties of this structures have been studied in (see [2] to [7], [11], [12], [15], [16]). We introduce the natural almost complex and product structures $\widetilde{J}$ and $\widetilde{Q}$ which depend only on $g$ and preserve the property of homogeneity. Then we get almost anti-Hermitian structure $(\tilde{g}_2, \widetilde{J})$ and almost product structure $(\tilde{g}_2, \widetilde{Q})$. By considering twin tensor of $\tilde{g}_2$, we construct almost para-Hermitian and Hermitian structures on $\widetilde{TM}$.

Let $M$ be a smooth manifold, $TM$ its tangent bundle and $\chi(M)$ the algebra of vector fields on $M$. A $K$-structure on $M$ is a fields of endomorphisms $K$ on $TM$ such that $K^2 = \varepsilon I$, where $\varepsilon = \pm 1$. Thus $\varepsilon = 1$ corresponds to an *almost product structure*, while $\varepsilon = -1$ provides an *almost complex structure*.

A $K$-structure is *integrable* if and only if there exists an a linear torsionless connection on $M$ such that $\nabla K = 0$, or equivalently the Nijenhuis tensor $N_K$ vanishes. In this case, $\nabla$ is called *almost complex (product) connection* if $K$ be an almost complex (product) structure.

If $g$ is a metric on $M$ such that $g(KX, KY) = \sigma g(X,Y)$, $\sigma = \pm 1$, for arbitrary vector fields $X$ and $Y$ on $M$, then we shall say that the metric $g$ is $K$-metric.

The definition above unifies the following four cases:

The case $\varepsilon = 1, \sigma = 1$ corresponds to the (pseudo-) Riemannian *almost product manifold* $(M, g, K)$, the case $\varepsilon = 1, \sigma = -1$ provides the *almost para-Hermitian manifold* $(M, g, K)$, the case $\varepsilon = -1, \sigma = 1$ is known as the *almost Hermitian manifold* $(M, g, K)$, and finally the case $\varepsilon = -1, \sigma = -1$ corresponds to the *almost anti-Hermitian* manifold $(M, g, K)$.

Let us introduce a $(0,2)$ tensor field $h$, the twin of $g$, by $h(X,Y) = g(KX,Y)$.
Then

$$h(X,Y) = \varepsilon\sigma h(Y,X), h(KX,KY) = \sigma h(X,Y).$$

Notice that for $\varepsilon\sigma = 1$, the twin tensor is a metric, while for $\varepsilon\sigma = -1$ the twin tensor is a 2-form.

Let $\psi$ be a (0,3) tensor fields defined by the formula

$$\psi(X,Y,Z) = g((\nabla_X K)Y, Z) \equiv (\nabla_X h)(Y,Z) \tag{1.1}$$

Obviously, if the tensor fields $\psi$ vanishes then $\nabla K = 0$ for a torsionless (Levi-Civita) connection and the Nijenhuis tensor $N_K$ is forced to vanish, too ([2]).

## 2. The Complete Lift

Let $\Gamma_{ij}^{\ k}$ be the coefficients of the Riemannian connection of $M$, then $N_j^{\ h} = \Gamma_{0j}^{\ h} = y^a \Gamma_{aj}^{\ h}(x)$ can be regarded as coefficients of the canonical nonlinear connection $N$ of $TM$, where $(x^h, y^h)$ are the induced coordinates in $TM$.

$N$ determines a horizontal distribution on $\widetilde{TM}$, which is supplementary to the vertical distribution $V$, such that, we have:

$$T_u\widetilde{TM} = N_u \oplus V_u, \quad \forall u \in \widetilde{TM}. \tag{2.1}$$

The adapted basis to $N$ and $V$ is given by $\{X_h, X_{\bar{h}}\}$ where

$$X_h = \frac{\partial}{\partial x^h} - y^a \Gamma_{ah}^{\ m} \frac{\partial}{\partial y^m}, \qquad X_{\bar{h}} = \frac{\partial}{\partial y^h} \tag{2.2}$$

and its dual basis is $\{dx^i, \delta y^i\}$ where

$$\delta y^i = dy^i + y^a \Gamma_{aj}^{\ i} dx^j. \tag{2.3}$$

The indices $a, b, ..., \bar{a}, \bar{b}, ...,$ run over the range $\{1, 2, ..., n\}$. The summation convention will be used in relation to this system of indices. By straightforward calculations, we have the following lemma.

**Lemma 1.** The Lie bracket of the adapted frame of $TM$ satisfies the following:

*1)* $[X_i, X_j] = y^a K_{jia}^{\ \ \ m} X_{\bar{m}}$,

2) $[X_i, X_{\bar{j}}] = \Gamma_{ji}^{\ m}$,

3) $[X_{\bar{i}}, X_{\bar{j}}] = 0$,

where $K_{jia}{}^{m}$ denote the components of the curvature tensor of $M$.
Let $(M, g)$ be a Riemannian space, $M$ being a real n-dimensional manifold and $(TM, \pi, M)$ its tangent bundle. On a domain $U \subset M$ of a local chart, $g$ has the components $g_{ij}(x)$, $(i, j, ... = 1, ..., n)$. Then on the domain of chart $\pi^{-1}(U) \subset TM$ we consider the functions $g_{ij}(x, y) = g_{ij}(x), \forall (x, y) \in \pi^{-1}(U)$ and put

$$\| y \| = \sqrt{g_{ij}(x) y^i y^j}. \tag{2.4}$$

Then, $\| y \|$ is globally defined on $TM$, differentiable on $\widetilde{TM}$ and continuous on the null section.

The complete lift of $g$ to $TM$ is defined by

$$g_2(x, y) = 2 g_{ij}(x) dx^i \delta y^j, \ \forall (x, y) \in \widetilde{TM}. \tag{2.5}$$

Then, $g_2$ is not 0-homogeneous on the fibers of $TM$.
Namely, for the homothety $h_t : (x, y) \to (x, ty)$ for all $t \in R^+$ we get

$$(g_2 \circ h_t)(x, y) = 2t g_{ij}(x) dx^i \delta y^j = t g_2(x, y) \neq g_2(x, y).$$

On $\widetilde{TM}$ we define an almost complex structure $J$ by

$$J(X_i) = -X_{\bar{i}}, \quad J(X_{\bar{i}}) = X_i, \quad i = 1, ..., n. \tag{2.6}$$

It is known that $(\widetilde{TM}, J, g_2)$ is an almost anti-Hermitian manifold. Moreover, the integrability of the almost complex structure $J$ implies that $(M, g)$ is locally flat. (see [7])
Also, we define almost product structure $Q$ on $\widetilde{TM}$ by

$$Q(X_i) = X_{\bar{i}}, \quad Q(X_{\bar{i}}) = X_i, \quad i = 1, ..., n. \tag{2.7}$$

Then, $(\widetilde{TM}, Q, g_2)$ is an almost product manifold. Also, the integrability of the almost product structure $Q$ implies that $(M, g)$ is locally flat.
The previous space, called "the geometrical model on $TM$ of the Riemannian space $(M, g)$", is important in the study of the geometry of initial Riamannian space $(M, g)$ ([6], [7]).

## 3. The 0-homogeneous lift of the Riemannian metric $g$

We can eliminate the inconvenience of the complete lift, introducing a new kind of lift to $TM$ of the Riemannian metric $g$. Then we obtain the Levi-Civita connection for this metric.
**Definition .** Let $\tilde{g}_2$ be a the tensor field on $\widetilde{TM}$ defined by

$$\tilde{g}_2(x, y) = \frac{2}{\| y \|} g_{ij}(x) dx^i \delta y^i \tag{3.1}$$

where $\| y \|$ was defined in (2.4).Then $\tilde{g}_2$ is called the 0-homogeneous lift of the Riemannian metric $g$ to $\widetilde{TM}$.
We get, evidently:
***Theorem 2.*** *The following properties hold:*
1. *The pair* $(\widetilde{TM}, \tilde{g}_2)$ *is a pseudo-Riemannian space, depending only on the metric* $g$.

2. $\tilde{g}_2$ *is 0-homogeneous on the fibers of the tangent bundle* $TM$.

In order to study the geometry of the pseudo-Riemannian space $(\widetilde{TM}, \tilde{g}_2)$ we can apply the theory of the $(h,v)$-Riemannian metric on $TM$ given in the books [6], [7] and [9]. Looking at the relation (2.5) and (3.1) we can assert:

***Theorem 3.*** *The lifts* $g_2$ *and* $\tilde{g}_2$ *coincide on the hyper unit tangent sphere* $g_{ij}(x_0)y^i y^j = 1$, *for every point* $x_0 \in M$.

Let $\overline{\nabla}$ be the Riemannian connection of $TM$ with coefficient $\overline{\Gamma}_{BC}{}^A$, that is:

$$\overline{\nabla}_{X_i} X_j = \overline{\Gamma}_{ji}{}^m X_m + \overline{\Gamma}_{ji}{}^{\bar{m}} X_{\bar{m}}, \qquad \overline{\nabla}_{X_i} X_{\bar{j}} = \overline{\Gamma}_{\bar{j}i}{}^m X_m + \overline{\Gamma}_{\bar{j}i}{}^{\bar{m}} X_{\bar{m}},$$
$$\overline{\nabla}_{X_{\bar{i}}} X_j = \overline{\Gamma}_{j\bar{i}}{}^m X_m + \overline{\Gamma}_{j\bar{i}}{}^{\bar{m}} X_{\bar{m}}, \qquad \overline{\nabla}_{X_{\bar{i}}} X_{\bar{j}} = \overline{\Gamma}_{\bar{j}\bar{i}}{}^m X_m + \overline{\Gamma}_{\bar{j}\bar{i}}{}^{\bar{m}} X_{\bar{m}} \tag{3.2}$$

Then, we have

$$\overline{\nabla}_{X_i} dx^h = -\overline{\Gamma}_{mi}{}^h dx^m - \overline{\Gamma}_{\bar{m}i}{}^h \delta y^m,$$
$$\overline{\nabla}_{X_i} \delta y^h = -\overline{\Gamma}_{mi}{}^{\bar{h}} dx^m - \overline{\Gamma}_{\bar{m}i}{}^{\bar{h}} \delta y^m,$$
$$\overline{\nabla}_{X_{\bar{i}}} dx^h = -\overline{\Gamma}_{m\bar{i}}{}^h dx^m - \overline{\Gamma}_{\bar{m}\bar{i}}{}^h \delta y^m, \tag{3.3}$$
$$\overline{\nabla}_{X_i} \delta y^h = -\overline{\Gamma}_{m\bar{i}}{}^{\bar{h}} dx^m - \overline{\Gamma}_{\bar{m}\bar{i}}{}^{\bar{h}} \delta y^m,$$

Since the torsion tensor $T(X,Y)$ of $\overline{\nabla}$ defined by $T(X,Y) = \overline{\nabla}_X Y - \overline{\nabla}_Y X - [X,Y]$ vanishes, we have the following relations by means of Lemma 1 and $(3.2)$.

$$(1)\ \overline{\Gamma}_{ji}{}^h = \overline{\Gamma}_{ij}{}^h \qquad (2)\ \overline{\Gamma}_{ji}^{\bar{h}} = \overline{\Gamma}_{ij}^{\bar{h}} + y^a K_{jia}^h$$
$$(3)\ \overline{\Gamma}_{\bar{j}i}{}^h = \overline{\Gamma}_{i\bar{j}}{}^h \qquad (4)\ \overline{\Gamma}_{\bar{j}i}{}^{\bar{h}} = \overline{\Gamma}_{i\bar{j}}{}^{\bar{h}} + \Gamma_{ji}{}^h \tag{3.4}$$
$$(5)\ \overline{\Gamma}_{\bar{j}\bar{i}}{}^h = \overline{\Gamma}_{\bar{i}\bar{j}}{}^h \qquad (6)\ \overline{\Gamma}_{\bar{j}\bar{i}}{}^{\bar{h}} = \overline{\Gamma}_{\bar{i}\bar{j}}{}^{\bar{h}}$$

Furthermore, we have the following lemma.

***Lemma 4.*** *The connection coefficients* $\overline{\Gamma}_{BC}{}^A$ *of* $\overline{\nabla}$ *of the complete metric* $\tilde{g}_2$ *satisfy the following relations:*

$$(1)\ \overline{\Gamma}_{ji}{}^h = \Gamma_{ji}{}^h \qquad (2)\ \overline{\Gamma}_{ji}{}^{\bar{h}} = y^a K_{aij}{}^h$$
$$(3)\ \overline{\Gamma}_{\bar{j}i}{}^h = \frac{1}{2\|y\|^2}(g_{ij}y^h - \delta_i^h y_j) \qquad (4)\ \overline{\Gamma}_{j\bar{i}}{}^h = \frac{1}{2\|y\|^2}(g_{ij}y^h - \delta_j^h y_i)$$
$$(5)\ \overline{\Gamma}_{\bar{j}i}{}^{\bar{h}} = \Gamma_{ji}{}^h \qquad (6)\ \overline{\Gamma}_{j\bar{i}}{}^{\bar{h}} = 0$$
$$(7)\ \overline{\Gamma}_{\bar{j}\bar{i}}{}^h = 0 \qquad (8)\ \overline{\Gamma}_{\bar{j}\bar{i}}{}^{\bar{h}} = -\frac{1}{2\|y\|^2}(\delta_i^h y_j + \delta_j^h y_i)$$

**Proof.** The condition compatibility $\overline{\nabla}$ is equivalent with following equations:

$$g_{ir}\overline{\Gamma}_{jm}{}^{\bar{r}} + g_{jr}\overline{\Gamma}_{im}{}^{\bar{r}} = 0 \tag{3.5}$$
$$g_{ir}(\overline{\Gamma}_{jm}{}^r - \overline{\Gamma}_{\bar{j}m}{}^{\bar{r}}) + g_{jr}(\Gamma_{im}{}^r - \overline{\Gamma}_{im}{}^r) = 0 \tag{3.6}$$
$$g_{ir}\overline{\Gamma}_{\bar{j}m}{}^r + g_{jr}\overline{\Gamma}_{\bar{i}m}{}^r = 0 \tag{3.7}$$

$$g_{ir}\bar{\Gamma}_{j\bar{m}}{}^{\bar{r}} + g_{jr}\bar{\Gamma}_{i\bar{m}}{}^{\bar{r}} = 0 \tag{3.8}$$

$$g_{ir}\bar{\Gamma}_{\bar{j}\bar{m}}{}^{\bar{r}} + g_{jr}\bar{\Gamma}_{i\bar{m}}{}^{r} + \frac{1}{\| y \|^2} g_{ij} y_m = 0 \tag{3.9}$$

$$g_{ir}\bar{\Gamma}_{\bar{j}\bar{m}}{}^{r} + g_{jr}\bar{\Gamma}_{\bar{i}\bar{m}}{}^{r} = 0 \tag{3.10}$$

From (3.10) we have $\bar{\Gamma}_{\bar{j}\bar{i}}{}^{h} = 0$, thus we get (7). From (3.4), (3.9) and (3.7), we have

$$g_{ir}\bar{\Gamma}_{\bar{j}m}{}^{r} = -g_{ir}\bar{\Gamma}^{r}_{\bar{i}m} = -g_{ir}\bar{\Gamma}^{r}_{m\bar{i}} = g_{mr}\bar{\Gamma}^{\bar{r}}_{\bar{i}\bar{j}} + \frac{g_{mj} y_i}{\| y \|^2}$$

$$= -g_{ir}\bar{\Gamma}^{r}_{m\bar{j}} - \frac{g_{im} y_j}{\| y \|^2} + \frac{g_{mj} y_i}{\| y \|^2}$$

$$= -g_{ir}\bar{\Gamma}^{r}_{\bar{j}m} - \frac{g_{im} y_j}{\| y \|^2} + \frac{g_{mj} y_i}{\| y \|^2}$$

Thus we get (3). From (3) and (3.4), we have (4).
From (3.9), (4) and (3.4), we have

$$g_{ir}\bar{\Gamma}_{\bar{j}\bar{m}}{}^{\bar{r}} + \frac{1}{2\| y \|^2} g_{im} y_j - \frac{1}{2\| y \|^2} g_{ji} y_m + \frac{g_{ij} y_m}{\| y \|^2} = 0,$$

then we obtain (8).
From (3.4) and (3.5) we have

$$g_{ir}\bar{\Gamma}_{jm}{}^{\bar{r}} = -g_{jr}\bar{\Gamma}_{im}{}^{\bar{r}} = -g_{jr}(\bar{\Gamma}_{mi}{}^{\bar{r}} + y^a K_{ima}{}^{r}) = g_{mr}\bar{\Gamma}_{ji}{}^{\bar{r}} - y^a K_{imaj} = g_{mr}\bar{\Gamma}_{ij}{}^{\bar{r}} + y^a (K_{jiam} - K_{imaj})$$

$$= -g_{ir}\bar{\Gamma}_{mj}{}^{\bar{r}} + y^a (K_{jiam} - K_{imaj}) = -g_{ir}(\bar{\Gamma}_{jm}{}^{\bar{r}} + y^a K_{mja}{}^{r}) + y^a (K_{jiam} - K_{imaj}),$$

thus we get (2).
From (3.4), (3.6) and (3.8), we have

$$g_{ir}\bar{\Gamma}_{j\bar{m}}{}^{\bar{r}} = -g_{jr}\bar{\Gamma}_{i\bar{m}}{}^{\bar{r}} = -g_{jr}(\bar{\Gamma}_{\bar{m}i}{}^{\bar{r}} - \Gamma_{mi}{}^{r}) = g_{jr}(\Gamma_{mi}{}^{r} - \bar{\Gamma}_{\bar{m}i}{}^{\bar{r}})$$

$$= -g_{mr}(\Gamma_{ji}{}^{r} - \bar{\Gamma}_{ji}{}^{r}) = -g_{mr}(\Gamma_{ij}{}^{r} - \bar{\Gamma}_{ij}{}^{r}) = g_{ir}(\Gamma_{mj}{}^{r} - \bar{\Gamma}_{\bar{m}j}{}^{\bar{r}})$$

$$= g_{ir}\Gamma_{mj}{}^{r} - g_{ir}\bar{\Gamma}_{\bar{m}j}{}^{\bar{r}} = g_{ir}\Gamma_{mj}{}^{r} - g_{ir}(\bar{\Gamma}_{j\bar{m}}{}^{\bar{r}} + \Gamma_{mj}{}^{r})$$

thus we obtain (5) and (6). From (3.6) and (5), we have (1).

## 4. The almost anti-Hermitian structure $(\tilde{g}_2, \tilde{J})$

The almost complex structure $J$ defined in (2.6) has not the property of homogeneity. The $F(\widetilde{TM})$-linear mapping $J : \chi(\widetilde{TM}) \to \chi(\widetilde{TM})$, applies the 1-homogeneous vector fields $X_i$ into 0-homogeneous vector fields $X_{\bar{i}}$ $(i = 1, \dots, n)$. Therefore, we consider the $F(\widetilde{TM})$-linear mapping $\tilde{J} : \chi(\widetilde{TM}) \to \chi(\widetilde{TM})$, given on the adapted basis by

$$\tilde{J}(X_i) = -\| y \| X_{\bar{i}}, \qquad \tilde{J}(X_{\bar{i}}) = \frac{1}{\| y \|} X_i, \ (i = 1, \dots, n). \tag{4.1}$$

Obviously, $\tilde{J}$ is a tensor field of type (1,1) on $\widetilde{TM}$, that is homogeneous on the fibers of $TM$.

***Theorem 5.*** $(\widetilde{TM}, \tilde{g}_2, \tilde{J})$ *is an almost anti-Hermitian manifold.*

**Proof**. It follows easily that

$$\tilde{g}_2(JX_i, JX_j) = -\tilde{g}_2(X_i, X_j), \quad \tilde{g}_2(JX_{\bar{i}}, JX_{\bar{j}}) = -\tilde{g}_2(X_{\bar{i}}, X_{\bar{j}})$$

$$\tilde{g}_2(JX_{\bar{i}}, JX_j) = -\tilde{g}_2(X_{\bar{i}}, X_j).$$

Hence

$$\tilde{g}_2(\tilde{J}X, \tilde{J}Y) = -\tilde{g}_2(X, Y), \ \forall X, Y \in \chi(\widetilde{TM}).$$

**Proposition 6.** *The Nijenhuis tensor field of the almost complex structure* $\tilde{J}$ *on* $\widetilde{TM}$ *is give by*

$$\begin{cases} N_{\tilde{J}}(X_i, X_j) = (y_i\delta_j^s - y_j\delta_i^s - y^a K_{jia}{}^s)X_{\bar{s}}, \\ N_{\tilde{J}}(X_i, X_{\bar{j}}) = \dfrac{1}{\|y\|^2}(y_i\delta_j^s - y_j\delta_i^s - y^a K_{jia}{}^s)X_s, \\ N_{\tilde{J}}(X_{\bar{i}}, X_{\bar{j}}) = \dfrac{1}{\|y\|^2}(y_j\delta_i^s - y_i\delta_j^s + y^a K_{jia}{}^s)X_{\bar{s}}. \end{cases}$$

**Proof**. Recall that the Nijenhuis tensor field $N_{\tilde{J}}$ defined by $\tilde{J}$ is given by

$$N_{\tilde{J}}(X, Y) = [\tilde{J}X, \tilde{J}Y] - \tilde{J}[\tilde{J}X, Y] - \tilde{J}[X, \tilde{J}Y] - [X, Y], \qquad \forall X, Y \in \chi(\widetilde{TM}).$$

Replacing the basis $(X_i, X_{\bar{i}})$ in the above formula and using following relation:

$$X_i(\|y\|) = 0, \quad X_{\bar{i}}(\|y\|) = \frac{y_i}{\|y\|}$$

We get the proof.

**Theorem 7.** *The almost complex structure* $\tilde{J}$ *is a complex structure on* $\widetilde{TM}$ *if and only if the Riemannian space* $(M, g)$ *is of constant sectional curvature 1.*

**Proof**. From the condition $N_{\tilde{J}} = 0$, one obtains:

$$\{K_{jia}{}^s - (g_{ia}\delta_j^s - g_{ja}\delta_i^s)\}y^a = 0.$$

Differentiating with respect to $y^h$, taking $y^a = 0 \ \ \forall a \in \{1, \ldots, n\}$, it follows that the curvature tensor field of $\nabla$ has the expression

$$K_{jih}{}^s = g_{ih}\delta_j^s - g_{jh}\delta_i^s \tag{4.2}$$

Using by the Schur theorem (in the case where $M$ is connected and $\dim M \geq 3$) it follows that $(M, g)$ has the constant sectional curvature 1.

**Corollary 8.** $(\widetilde{TM}, \tilde{g}_2, \tilde{J})$ *is an anti-Hermitian manifold if and only if the space* $(M, g)$ *is of constant sectional curvature 1.*

From $(4.2)$ we have

$$R_{ij} = (n-1)g_{ij}, \qquad (n > 1) \tag{4.3}$$

where $R_{rk}$ is the Ricci tensor and

$$S = n(n-1).$$

where $S$ is the scalar tensor.

**Corollary 9.** *If the structure* $(\tilde{g}_2, \widetilde{J})$ *is a Hermitian structure on* $\widetilde{TM}$ *then* $(M, g)$ *is an Einstein space with positive scalar curvature.*

Since $R_{ij} = R_{ji}$ then from $(4.3)$ we get:

**Corollary 10.** *If the almost complex structure* $\widetilde{J}$ *is a complex structure then* $(M, R_{ij}(x))$ *is a Riemannian space.*

## 5. The almost product structure $(\tilde{g}_2, \widetilde{Q})$

The almost product structure $Q$ defined in $(2.7)$ has not the property of homogeneity. The $F(\widetilde{TM})$-linear mapping $Q : \chi(\widetilde{TM}) \to \chi(\widetilde{TM})$, applies the 1-homogeneous vector fields $X_i$ into 0-homogeneous vector fields $X_{\bar{i}}\,(i = 1, ..., n)$. Therefore, we consider the $F(\widetilde{TM})$-linear mapping $\widetilde{Q} : \chi(\widetilde{TM}) \to \chi(\widetilde{TM})$, given on the adapted basis by

$$\widetilde{Q}(X_i) = \| y \| X_{\bar{i}}, \quad \widetilde{Q}(X_{\bar{i}}) = \frac{1}{\| y \|} X_i, \ (i = 1, ..., n). \tag{5.1}$$

Obviously, $\tilde{Q}$ is a tensor field of type $(1,1)$ on $\widetilde{TM}$, that is homogeneous on the fibers of $TM$. It is not difficult to prove:

**Theorem 11.** $(\widetilde{TM}, \tilde{g}_2, \widetilde{Q})$ *is an almost product manifold.*

In order to find conditions that $\widetilde{Q}$ be a product structure, we have to put zero for the Nijenhuis tensor field of $\widetilde{Q}$,

$$N_{\tilde{Q}}(X, Y) = [\tilde{Q}X, \tilde{Q}Y] - \tilde{Q}[\tilde{Q}X, Y] - \tilde{Q}[X, \tilde{Q}Y] + [X, Y], \qquad \forall X, Y \in \chi(\widetilde{TM}).$$

**Theorem 12.** $(\widetilde{TM}, \tilde{g}_2, \tilde{Q})$ *is a product manifold if and only if the space* $(M, g)$ *is of constant sectional curvature -1.*

**Proof.** Similar to proposition 6 and theorem 7, by putting $N_{\tilde{Q}} = 0$ we get

$$K_{jia}{}^{s} = -(g_{ia}\delta^s_j - g_{ja}\delta^s_i). \tag{5.2}$$

Therefore, using by the Schur theorem, it follows that $(M, g)$ has the constant sectional curvature -1.

**Theorem 13.** *If the structure* $(\tilde{g}_2, \widetilde{Q})$ *is a product structure on* $\widetilde{TM}$ *then* $(M, g)$ *is an Einstein space with negative scalar curvature.*

**Proof.** From $(5.2)$ we have $R_{ij} = (1-n)g_{ij}, S = n(1-n)$ for $n > 1$.

Since $R_{ij} = R_{ji}$ then we get:

**Corollary 14.** *If the almost product structure* $\widetilde{Q}$ *is a product structure then* $(M, R_{ij}(x))$ *is a Riemannian space.*

## 6. Almost Hermitian and para-Hermitian structures on $\widetilde{TM}$

In this section, we get twin tensor of metric $\tilde{g}_2$ and by using it introduce Almost Hermitian and para-Hermitian structures on $\widetilde{TM}$. Then we show that these structures are not Kahlerian or para-Kahlerian.

**Lemma15.** *The twin tensor of structure* $(\tilde{g}_2, \tilde{J})$ *is a metric that is given by*

$$h_{\tilde{J}} = -2g_{ij}dx^i dx^j + \frac{2g_{ij}}{\|y\|^2}\delta y^i \delta y^j$$

**Proof.** From relation $h_{\tilde{J}}(X,Y) = \tilde{g}_2(\tilde{J}X,Y)$ we have:

$$h_{\tilde{J}}(X_i, X_j) = -\|y\|\,\tilde{g}_2(X_{\bar{i}}, X_j) = -2g_{ij},\ h_{\tilde{J}}(X_{\bar{i}}, X_{\bar{j}}) = \frac{1}{\|y\|}\tilde{g}_2(X_i, X_{\bar{j}}) = \frac{2}{\|y\|^2}g_{ij}$$

$$h_{\tilde{J}}(X_i, X_{\bar{j}}) = -\|y\|\,\tilde{g}_2(X_{\bar{i}}, X_{\bar{j}}) = 0.$$

**Theorem 16.** $(\widetilde{TM}, h_{\tilde{J}}, \tilde{Q})$ *is an almost para-Hermitian manifold.*

**Proof.** Straightforward computations, we obtain

$$h_{\tilde{J}}(\tilde{Q}X_i, \tilde{Q}X_j) = 2g_{ij} = -h_{\tilde{J}}(X_i, X_j),\ h_{\tilde{J}}(\tilde{Q}X_{\bar{i}}, \tilde{Q}X_{\bar{j}}) = -\frac{2}{\|y\|^2}g_{ij} = -h_{\tilde{J}}(X_{\bar{i}}, X_{\bar{j}})$$

$$h_{\tilde{J}}(\tilde{Q}X_{\bar{i}}, \tilde{Q}X_j) = 0 = -h_{\tilde{J}}(X_{\bar{i}}, X_j)$$

Therefore

$$h_{\tilde{J}}(\tilde{Q}X, \tilde{Q}Y) = -h_{\tilde{J}}(X,Y)$$

By definition $\Omega_{\tilde{Q}}^{h_{\tilde{J}}}(X,Y) = h_{\tilde{J}}(\tilde{Q}X,Y)$, the associated almost simplectic structure $\Omega^{h_{\tilde{J}}}$ is given in adapted basis by

$$\Omega^{h_{\tilde{J}}} = \frac{4}{\|y\|}g_{ij}(x)dx^i \wedge \delta y^j.$$

**Theorem17.** *The space* $(\widetilde{TM}, h_{\tilde{J}}, \tilde{Q})$ *cannot be an almost para- Kählerian manifold.*

**Proof.** since, $d(\frac{1}{\|y\|}) = -\frac{1}{\|y\|^2}d\|y\|$ and $d(g_{ij}dx^i \wedge \delta y^j) = 0$ then, the exterior differential of $\Omega^{h_{\tilde{J}}}$ satisfies the equation:

$$d\Omega^{h_{\tilde{J}}} = -\frac{4}{\|y\|}d\|y\| \wedge \Omega^{h_{\tilde{J}}}.$$

It follows, easily that $d\Omega_{\tilde{J}}^{h_{\tilde{Q}}} \neq 0$ on $\widetilde{TM}$, i.e., $\Omega_{\tilde{J}}^{h_{\tilde{Q}}}$ is not closed.

From theorem 12,16, we have:

**Theorem18.** $(\widetilde{TM}, h_{\tilde{J}}, \tilde{Q})$ *is a para-Hermitian manifold if and only if the space* $(M,g)$ *is of constant sectional curvature -1.*

**Lemma19.** *The Levi-Civita connection coefficients* $\overline{\nabla}^{h_{\tilde{J}}}$ *of* $h_{\tilde{J}}$ *satisfy the following relations:*

(1) $\bar{\Gamma}_{ji}{}^{h} = \Gamma_{ji}{}^{h}$, (2) $\bar{\Gamma}_{ji}{}^{\bar{h}} = \frac{1}{2} y^a K_{jia}{}^{h}$,

(3) $\bar{\Gamma}_{\bar{j}i}{}^{h} = -\frac{1}{2\|y\|^2} y^a K_{aji}{}^{h}$, (4) $\bar{\Gamma}_{j\bar{i}}{}^{h} = -\frac{1}{2\|y\|^2} y^a K_{aij}{}^{h}$,

(5) $\bar{\Gamma}_{\bar{j}i}{}^{\bar{h}} = \Gamma_{ji}{}^{h}$, (6) $\bar{\Gamma}_{j\bar{i}}{}^{\bar{h}} = 0$,

(7) $\bar{\Gamma}_{\bar{j}\bar{i}}{}^{h} = 0$, (8) $\bar{\Gamma}_{\bar{j}\bar{i}}{}^{\bar{h}} = \frac{1}{\|y\|^2}(g_{ji} y^h - \delta_j^h y_i - \delta_i^h y_j)$.

**Theorem 20.** $\bar{\nabla}^{h_{\tilde{J}}}$ *is an almost complex connection.*
**Proof.** From $(1.1)$ we have

$$\tilde{g}_2((\bar{\nabla}_X^{h_{\tilde{J}}} \tilde{J})Y, Z) \equiv (\bar{\nabla}_X^{h_{\tilde{J}}} h_{\tilde{J}})(Y, Z)$$

Since $\bar{\nabla}^{h_{\tilde{J}}}$ is Levi-Civita connection for $h_{\tilde{J}}$ then

$$\tilde{g}_2((\bar{\nabla}_X^{h_{\tilde{J}}} \tilde{J})Y, Z) = 0$$

i.e. $\nabla_X \tilde{J} = 0$.

Similarly previous case, the twin tensor of structure $(\tilde{g}_2, \tilde{Q})$ is a metric that is

$$h_{\tilde{Q}} = 2 g_{ij} dx^i dx^j + \frac{2 g_{ij}}{\|y\|^2} \delta y^i \delta y^j$$

Obviously, $h_{\tilde{Q}}$ is 0-homogeneous on the fibers of TM.

**Theorem 21.**
*1.* $(\widetilde{TM}, h_{\tilde{Q}}, \tilde{J})$ *is an almost Hermitian structure on* $\widetilde{TM}$ *.*

*2. The associated almost simplectic structure* $\Omega_{\tilde{J}}^{h_{\tilde{Q}}}$ *is given in adapted basis by*

$$\Omega_{\tilde{J}}^{h_{\tilde{Q}}} = \frac{4}{\|y\|} g_{ij}(x) \delta y^i \wedge dx^j$$

**Theorem 22.** *The space* $(\widetilde{TM}, h_{\tilde{Q}}, \tilde{J})$ *cannot be an almost Kählerian manifold.*

From corollary 8 and theorem 21 we obtain following theorem.

**Theorem 23.** $(\widetilde{TM}, h_{\tilde{Q}}, \tilde{J})$ *is a Hermitian manifold if and only if the space* $(M, g)$ *is of constant sectional curvature 1.*

**Lemma 24.** *The Levi-Civita connection coefficients* $\bar{\nabla}^{h_{\tilde{Q}}}$ *of* $h_{\tilde{Q}}$ *satisfy the following relations:*

(1) $\bar{\Gamma}_{ji}{}^{h} = \Gamma_{ji}{}^{h}$, (2) $\bar{\Gamma}_{ji}{}^{\bar{h}} = \frac{1}{2} y^a K_{jia}{}^{h}$,

(3) $\bar{\Gamma}_{\bar{j}i}{}^{h} = \frac{1}{2\|y\|^2} y^a K_{aji}{}^{h}$, (4) $\bar{\Gamma}_{j\bar{i}}{}^{h} = \frac{1}{2\|y\|^2} y^a K_{aij}{}^{h}$,

(5) $\bar{\Gamma}_{\bar{j}i}{}^{\bar{h}} = \Gamma_{ji}{}^{h}$, (6) $\bar{\Gamma}_{j\bar{i}}{}^{\bar{h}} = 0$,

(7) $\bar{\Gamma}_{\bar{j}\bar{i}}{}^{h} = 0$, (8) $\bar{\Gamma}_{\bar{j}\bar{i}}{}^{\bar{h}} = \frac{1}{\|y\|^2}(g_{ji} y^h - \delta_j^h y_i - \delta_i^h y_j)$.

**Theorem 25.** $\bar{\nabla}^{h_{\tilde{Q}}}$ *is an almost product connection.*

**Acknowledgement:** The authors would like to thanks Professor Radu Miron for advising to work on this field.

Department of Mathematics and Computer Science
AmirKabir University.
Tehran.Iran.
E-mail address: e_peyghan@aut.ac.ir
E-mail address: arazavi@aut.ac.ir

Faculty of Science of Tarbiatmodares University.
Tehran.Iran.
E-mail address: abasheydari@gmail.com